\documentclass[12pt]{amsart}
\usepackage{fullpage}
\usepackage{amsmath,amssymb,amscd}

\newcommand{\sA}{\mathcal{A}}
\newcommand{\bC}{\mathbb{C}}
\newcommand{\bQ}{\mathbb{Q}}
\newcommand{\bZ}{\mathbb{Z}}
\newcommand{\indlim}{\underrightarrow{ \mbox{\rm lim\vphantom{g}}}\,}
\newcommand{\indalglim}{\underrightarrow{\mbox{\rm alg\,lim}}\,}
\newcommand{\hr}{\mbox{hr}\,}
\newcommand{\rad}{\mbox{rad}\,}
\newcommand{\cstar}{\textrm{C}$^*$}

\def\symbup#1{\Big\uparrow\rlap{$\vcenter{\hbox{$\scriptstyle #1$}}$}}
\def\symbupdown#1{\Big\updownarrow\rlap{$\vcenter{\hbox{$\scriptstyle #1$}}$}}

\def\symbse#1{\!\!\!\!\searrow\rlap{$\vcenter{\hbox{$\scriptstyle #1$}}$}}

\theoremstyle{plain}
\newtheorem{thm}{Theorem}
\numberwithin{thm}{section}
\newtheorem{cor}[thm]{Corollary}
\newtheorem{lma}[thm]{Lemma}

\theoremstyle{definition}
\newtheorem{example}[thm]{Example}
\newtheorem{defn}[thm]{Definition}

\begin{document}
\bibliographystyle{amsplain}

\author[A.P. Donsig]{Allan P. Donsig}
\address{Dept. of Mathematics \& Statistics\\
	 University of Nebraska---Lincoln\\
	 Lincoln, NE,  68588-0323\\
	 U.S.A.}
\email{adonsig@math.unl.edu}
\author[S.C. Power]{S.C. Power}
\address{Dept. of Mathematics \& Statistics\\
	 Lancaster University\\
	 Lancaster, LA1 4YF\\
	 United Kingdom}
\email{s.power@lancaster.ac.uk}
\title{The classification of limits of $2n$-cycle algebras}
\thanks{Partial support was provided by
	an EPSRC grant and a UNL Summer Research Fellowship.}
\thanks{17 October 1998}
\date{}
\maketitle

\section{Introduction}
In this paper, we obtain a complete classification
of the locally finite algebras $\sA_0 = \indalglim A_k$ and
the operator algebras $\sA = \indlim A_k$
associated with towers
\[ A_1 \subseteq A_2 \subseteq A_3 \cdots \]
consisting of $2n$-cycle algebras, where $n \ge 3$,
with the inclusions of rigid type.
The complete isomorphism invariant is essentially
the triple \[ (K_0(\sA), H_1(\sA), \Sigma(\sA)) \]
where $K_0(\sA)$ is viewed as a scaled ordered group,
$H_1(\sA)$ is a partial isometry homology group
and $\Sigma(\sA) \subseteq K_0(\sA) \oplus H_1(\sA)$
is the $2n$-cycle joint scale.

Recall that a cycle algebra is a finite dimensional digraph algebra,
or incidence algebra, whose reduced digraph is a cycle.
For example, the basic $6$-cycle algebra in $M_6$ has the form
\[ \begin{bmatrix} * & * & & & & * \\ & * \\ & * & * & * \\ & & & * \\
	& & & * & * & * \\ & & & & & * \end{bmatrix} \]
and its digraph is a hexagon with alternating edge directions
and a loop edge at each vertex.
These algebras are of interest in that they are
the simplest family of finite dimensional
complex algebras with non-zero homology.
We focus on the
natural embeddings between cycle algebras known as rigid embeddings
and which have the property that they are determined by $K_0$ and $H_1$.
More precisely, a rigid embedding $\phi$ is determined
up to inner conjugacy by
$K_0 \phi \oplus H_1 \phi$, where $H_1 \phi$ is the
induced homomorphism between the first integral simplicial
homology groups of the simplicial complexes affiliated
to the digraphs of the algebras.
For cycle algebras, the map $H_1\phi$ is simply a group homomorphism
$\bZ \to \bZ$. This $K_0H_1$ uniqueness is in analogy with the fact that
embeddings between finite dimensional semisimple complex algebras
are determined up to inner conjugacy by $K_0$.

Limit homology groups for direct limits of digraph algebras were
introduced in~\cite{DavPow91} and an intrinsic formulation as stable
partial-isometry homology groups was given in~\cite{Pow96a}.
The classification of cycle algebras was first considered
in~\cite[Chapter~11]{Power92}, in the restricted context of
direct limits of 4-cycle algebras with
homologically limited embeddings.
These results were then extended to limits of 4-cycle algebras with
general rigid embeddings in~\cite{DonPow97}, and the classification
of the operator algebras up to {\sl regular} isomorphism
was obtained in terms of the triple $(K_0(\sA), H_1(\sA), \Sigma(\sA))$.
In~\cite{PowPRc}, the second author showed that this classification was, 
in fact, up to star extendible isomorphism for all locally finite algebras 
and for the operator algebras arising from direct systems of
4-cycle algebras in the somewhat more tractable odd case.
The current  paper is concerned with higher cycles and,
building on the ideas of~\cite{DonPow97}  
and~\cite{PowPRc}, we classify, 
up to star-extendible isomorphism, both the locally finite algebras 
and the operator algebras arising from direct limits of $2n$-cycle algebras, 
$n \ge 3$, with rigid embeddings.

Cohomology has been considered in the context of
non-selfadjoint algebras for many years;
recent references include~\cite{GilPopSmi97,GilSmi94}.
We remark that in addition to homological augmentations of $K_0$
invariants one can also consider scaled Grothendieck group
invariants for regular systems of digraph
algebras and their operator algebras. This topic is developed
in \cite{Pow97}.

It is fortuitous that, in one important respect, the analysis of
$2n$-cycle algebras for $n \ge 3$ is much simpler than that of $4$-cycle
algebras: star-extendible isomorphisms between algebraic direct limits
are necessarily induced by a commuting diagram of regular linking maps
between the given towers.
A \textsl{regular} morphism is a direct sum of multiplicity one embeddings.
The existence of irregular morphisms between regular limit
algebras was first observed in~\cite{DonPow96}.
On the other hand, the $K_0$-groups of towers of $2n$-cycle
algebras are much less readily identifiable, except in some
cases of independent interest, most notably, when $\sA \cap \sA^*$
is a direct sum of matroid algebras.

\textbf{Organisation}.
In the next section, we study $2n$-cycle algebras and the rigid embeddings
between them, culminating in the theorem that locally finite algebras
are star-extendibly isomorphic if and only if there is an intertwining
rigid diagram.
In the third section, we extend this theorem from locally finite algebras
to norm-closed direct limits.
The last section introduces the classifying invariants and obtains
the main classification theorem.
Finally we illustrate this with an application to a family of 
limit algebras for which the only variations are ones of homology.

\section{Towers of cycle algebras}

We first set out notation and recall some terminology.

A \textsl{digraph} is a finite directed graph with no multiple directed
edges of the same orientation.
To a transitive, reflexive digraph $G$ on the vertices $v_1,\ldots, v_n$,
we associate the complex algebra $A(G)$ spanned by those
standard matrix units $e_{i,j} \in M_n(\bC)$
associated with edges in $G$ from $v_j$ to $v_i$.
These subalgebras of $M_n(\bC)$ are precisely those which
contain the diagonal algebra.

For convenience we restrict attention henceforth to those morphisms
between digraph algebras $A(G_1) \to A(G_2)$ which are star-extendible,
in the sense of being restrictions of star-algebra homomorphisms
$C^*(A(G_1)) \to C^*(A(G_2))$.

A \textit{$2n$-cycle digraph algebra} is a digraph algebra $A(G)$
for which the reduced digraph of $G$, call it $G_r$, is isomorphic
to the connected graph on $2n$ vertices with the $2n$ edges
$(i,i)$ for $i=1,\ldots,2n$ and $2n$ edges between successive
vertices, with alternating orientations.
We may assume that the vertex labeled $1$ is a range vertex
rather than a source vertex.
Denote this directed graph by $D_{2n}$.

The elements of a $2n$-cycle algebra thus have a block matrix
staircase form:

\[ a = \begin{bmatrix} a_{11} & a_{12} & & & & & & a_{1,2n} \\
	& a_{22} \\
	& a_{32} & a_{33} & a_{34} \\
	& & & a_{44} \\
	& & & a_{54} & a_{55} & a_{56} \\
	& & & & & \ddots \\
	& & & & & & \ddots  & a_{2n-1,2n} \\
        & & & & & & & a_{2n,2n}
	\end{bmatrix}. \]

A \textsl{partial isometry} is an element $u$ for which $u^*u$
(and hence $u u^*$) is a projection (a selfadjoint idempotent).
A convenient consequence of the star-extendibility of an
embedding $\phi : A_1 \to A_2$ between digraph algebras is
that, for each partial isometry $v \in A_1$ with $vv^*$ and
$v^*v$ in $A_1$, the image $\phi(v)$ is a partial isometry
with its initial and final projections in $A_2$.
In $2n$-cycle algebras with $n \ge 3$, such partial isometries
have a particularly clear form.
We remark that the next lemma does not hold for $4$-cycle algebras,
which complicates considerably the analysis of their limit algebras.

\begin{lma}
\label{lma:PIs}
Suppose $a = \bigl( a_{ij} \bigr)$ is a block matrix
in a $2n$-cycle algebra $A$, $n \ge 3$.
If $a$ is a partial isometry with initial and final projections
in $A$, then each entry $a_{ij}$ is a partial isometry.
\end{lma}

\begin{proof}
The $(1,1)$ block entry of $a^*a$ is $a_{11}^* a_{11}$.
Thus, by hypothesis $a_{11}^* a_{11}$ is a projection and so
$a_{11}$ is a partial isometry.
By symmetry, $a_{kk}$ is a partial isometry for all $k$.
It follows that the block matrix
\[ b = \begin{bmatrix} a_{12} & & & & a_{1,2n} \\
	a_{32} & a_{34} \\
	& a_{54} & a_{56} \\
	& & \ddots & \ddots \\
	& & & a_{2n-1,2n-2} & a_{2n-1,2n}
	\end{bmatrix} \]
is a partial isometry, and that its initial and final projections
are block diagonal.
Since the $(1,1)$ entry of $bb^*$ is zero, it follows that
$a_{12}a_{32}^*=0$ and hence that $a_{12}$ and $a_{32}$ have
orthogonal initial projections.
Similarly, the $(1,2n)$ entry of $b^*b$ being zero implies that
$a_{12}$ and $a_{1,2n}$ have orthogonal final projections.
Since $b$ is a partial isometry, this double orthogonality
forces $a_{12}$ to be a partial isometry.
By symmetry, each $a_{i\,j}$ is a partial isometry,
for all $i$ and $j$.
\end{proof}

In a general digraph algebra $A$, we refer to a partial
isometry $v$ as a \textsl{regular partial isometry}
if $pvq$ is a partial isometry for each pair of
central projections
$p,q$ in $A \cap A^*$.
In particular, such a partial isometry can be written as
a sum of rank one partial isometries in $A$.
Note that a regular embedding sends regular partial
isometries to regular partial isometries.

We call an embedding \textsl{locally regular} if it maps
regular partial isometries to regular partial isometries.
This is a strictly weaker property, as we now show.

\begin{example}
\label{example:notlocreg}
Consider the upper triangular realisation of $A(D_4)$ spanned
by the diagonal matrix units and the
4-cycle $\{e_{13},e_{14},e_{24},e_{23}\}$.
Define $\phi$ from $A(D_4)$ to $A(D_4) \otimes M_4$ by
mapping this rank-one 4-cycle in $A(D_4)$ to
the 4-cycle $\{ v_1, v_2, v_3, v_4 \}$ where, as
$2 \times 2$ block matrices, the $v_i$ have the following
upper right blocks and are otherwise zero:
\[ \left[ \begin{array}{cc|cc|cc|cc}
	1 & 0 & & & 0 & 0 \\ 0 & 0 & & & 1 & 0 \\ \hline
	\phantom{0} & \phantom{0} & \phantom{0} & \phantom{0} & \phantom{0} & 
	\phantom{0} & \phantom{0} & \phantom{0} \\
	\phantom{0} & \phantom{0} & \phantom{0} & \phantom{0} & \phantom{0} & 
	\phantom{0} & \phantom{0} & \phantom{0} \\ \hline
	0 & 1 & & & 0 & 0 \\ 0 & 0 & & & 0 & 1 \\ \hline
	\phantom{0} & \phantom{0} & \phantom{0} & \phantom{0} & \phantom{0} & 
	\phantom{0} & \phantom{0} & \phantom{0} \\
	\phantom{0} & \phantom{0} & \phantom{0} & \phantom{0} & \phantom{0} & 
	\phantom{0} & \phantom{0} & \phantom{0} \\ \end{array} \right],
   \frac{1}{\sqrt{2}} \left[ \begin{array}{cc|cc|cc|cc}
	& & 1 & 0 & & & 1 & 0 \\ & & 1 & 0 & & & -1 & 0 \\ \hline
	\phantom{0} & \phantom{0} & \phantom{0} & \phantom{0} & \phantom{0} & 
	\phantom{0} & \phantom{0} & \phantom{0} \\ 	
	\phantom{0} & \phantom{0} & \phantom{0} & \phantom{0} & \phantom{0} & 
	\phantom{0} & \phantom{0} & \phantom{0} \\ \hline
	& & 0 & 1 & & & 0 & 1 \\ & & 0 & 1 & & & 0 & -1 \\ \hline
	\phantom{0} & \phantom{0} & \phantom{0} & \phantom{0} & \phantom{0} & 
	\phantom{0} & \phantom{0} & \phantom{0} \\
	\phantom{0} & \phantom{0} & \phantom{0} & \phantom{0} & \phantom{0} & 
	\phantom{0} & \phantom{0} & \phantom{0} \\ \end{array} \right], \]
and
\[  \frac{1}{\sqrt{2}} \left[ \begin{array}{cc|cc|cc|cc}
	\phantom{0} & \phantom{0} & \phantom{0} & \phantom{0} & \phantom{0} & 
	\phantom{0} & \phantom{0} & \phantom{0} \\
	\phantom{0} & \phantom{0} & \phantom{0} & \phantom{0} & \phantom{0} & 
	\phantom{0} & \phantom{0} & \phantom{0} \\ \hline
	& & 1 & 1 & & &0 &0 \\ & & 0 & 0 & & & 1 & 1 \\ \hline
	\phantom{0} & \phantom{0} & \phantom{0} & \phantom{0} & \phantom{0} & 
	\phantom{0} & \phantom{0} & \phantom{0} \\ 	
	\phantom{0} & \phantom{0} & \phantom{0} & \phantom{0} & \phantom{0} & 
	\phantom{0} & \phantom{0} & \phantom{0} \\ \hline
	& & 1 & -1 & & & 0 & 0 \\ & & 0 & 0 & & & 1 & -1 \\ \end{array}  
\right],
    \left[ \begin{array}{cc|cc|cc|cc}
	\phantom{0} & \phantom{0} & \phantom{0} & \phantom{0} & \phantom{0} & 
	\phantom{0} & \phantom{0} & \phantom{0} \\
	\phantom{0} & \phantom{0} & \phantom{0} & \phantom{0} & \phantom{0} & 
	\phantom{0} & \phantom{0} & \phantom{0} \\ \hline
	1 & 0 & & & 0 & 0 \\ 0 & 0 & & & 1 & 0 \\ \hline
	\phantom{0} & \phantom{0} & \phantom{0} & \phantom{0} & \phantom{0} & 
	\phantom{0} & \phantom{0} & \phantom{0} \\
	\phantom{0} & \phantom{0} & \phantom{0} & \phantom{0} & \phantom{0} & 
	\phantom{0} & \phantom{0} & \phantom{0} \\ \hline
	0 & 1 & & & 0 & 0 \\0 & 0 & & & 0 & 1 \\ \end{array} \right]. \]
Since the $v_i$ are regular partial isometries with orthogonal
initial and final projections, $\phi$ is locally regular and star-extendible.
However, the product $v_2v_1^*$ is easily seen to be not locally regular.
It follows that $\phi$ is not regular.
\end{example}

\subsection{Rigid Embeddings}

The digraph $D_{2n}$ has $2n$ automorphisms which induce
$2n$ automorphisms of $A(D_{2n})$, denoted
$\theta_1, \theta_2, \cdots, \theta_{2n}$.
For definiteness, we let $\theta_1$ be the identity, $\theta_2$ be the
reflection which fixes the vertex $1$, $\theta_3$ be the shift which
maps each vertex $k$ to $(k-2) \mod 2m$, $\theta_{2k-1}=\theta_3^{k-1}$,
for $1 \le k \le m$, and $\theta_{2k} = \theta_2 \theta_{2k-1}$, for
$2 \le k \le m$.

Suppose that $A$ is a $2n$-cycle algebra and $i : A(D_{2n}) \to A$
is a multiplicity one star-extendible embedding, which
is \textsl{proper}, in the sense that $i(e_{jj})$, $1 \le j \le 2n$,
are inequivalent projections.
Then the embeddings
$i \circ \theta_1, i \circ \theta_2, \ldots,
	i \circ \theta_{2n}$
are representatives for the $2n$ inner unitary equivalence
classes of the proper multiplicity one injections.
We say that a star-extendible embedding
$\phi : A(D_{2n}) \to A(D_{2n}) \otimes M_n(\bC)$
is \textsl{rigid} if it decomposes as direct sum of
multiplicity one proper embeddings.

In general, a \textsl{rigid} 
embedding between $2n$-cycle algebras, $\phi : A_1 \to A_2$,
is one for which $\phi \circ \eta$ is rigid 
whenever $\eta: A(D_{2n}) \to A_1$ is proper.

We now come to the crucial lemma in the analysis of rigid embeddings.
The case $n=2$ was proved in Lemma~3.2 of \cite{PowPRc}
using a similar, albeit simpler, argument.

\begin{lma}
\label{lma:rigidfactors}
Let $\phi : A_1 \to A_2$ and $\psi : A_2 \to A_3$ be
locally regular (star-extendible) embeddings between
$2n$-cycle algebras, where $n \ge 3$.
If $\psi \circ \phi$ is rigid, then $\phi$ and $\psi$
are rigid.
\end{lma}

\begin{proof}

Suppose for simplicity of notation, that $n=3$.
The following argument does not depend on the length of
the cycle, and so it suffices to prove the lemma.

\medskip
Let $A_1$, $A_2$, $A_3$ be 6-cycle algebras and let
$\phi : A_1 \to A_2$ and $\psi : A_2 \to A_3$ be locally
regular embeddings whose composition, $\psi \circ \phi$,
is rigid.

We may assume that $\phi$ and $\psi$ are unital and that
$A_1=A(D_6)$ as the general case follows readily from this.
It is convenient to change our block matrix form for 6-cycle
algebras from a staircase pattern to the following one:
\[ \begin{bmatrix} * & & & * & & * \\ & * & & * & * \\ & & * & & * & * \\
	& & & * \\ & & & & * \\ & & & & & * \end{bmatrix} \]
Denote the partial isometries of the usual 6-cycle in $A(D_6)$
by $E_1, E_2, \ldots, E_6$ where these have the positions
$$
\left[ \begin{array}{ccc} E_1 & & E_6 \\ E_2 & E_3 & \\
	& E_4 & E_5 \end{array} \right]
$$
in the off-diagonal block of $A(D_6)$.
Since the composition $\psi \circ \phi$ maps each $E_i$ to an element
supported in the off diagonal blocks of $A_3$, it follows that $\phi$
acts similarly.
Thus, with respect to the block structure of $A_2$, we can write
$$
\phi(E_i) = \left[ \begin{array}{ccc} a_i & & f_i \\ b_i & c_i & \\
	& d_i & e_i \end{array} \right]
$$
for $i=1,2,\ldots, 6$.
By the local regularity hypothesis, each of the entries
$a_i,\ldots, f_i$ is
a partial isometry.
Moreover, by star-extendibility, the initial projection of $a_i$ is
orthogonal
to that of $b_i$ and, trivially, to those of $c_i, \ldots, f_i$.

Similarly, if $F_1,\ldots, F_6$ is a rank-one six-cycle in $A_2$ then,
with respect to the block structure of $A_3$, we can write
$$
\psi(F_i) = \left[ \begin{array}{ccc} \alpha_i & & \lambda_i \\ \beta_i  
	& \gamma_i & \\ & \delta_i & \epsilon_i \end{array} \right]
$$
for $i=1,\ldots, 6$.
As with $\phi$, local regularity implies that $\alpha_i, \ldots,
\lambda_i$ are partial
isometries and star-extendibility implies that, for example, $\alpha_i$  
and $\beta_i$
have orthogonal initial projections.

Fix matrix unit systems $\{e_{i,j}\}, \{f_{i,j}\}$ for $A_1$ and $A_2$, 
and note that we can assume that $\phi$ (resp. $\psi$) maps
matrix units in the self-adjoint
algebra $A_1 \cap A_1^*$ (resp. $A_2 \cap A_2^*$) to sums of matrix
units in $A_2$ (resp. $A_3$) and also that the restrictions
$\phi | A_1 \cap A_1^*$ and $ \psi | A_2 \cap A_2^*$ are standard embeddings.
(This may be arranged by replacing $\phi, \psi$ by
inner conjugate maps.)
We may assume for definiteness that $F_1, \dots F_6$ are the matrix units that
appear in the top left position of their block.
Thus $a_1$ can be written as
\[
a_1 = \sum a^1_{i,j}f_{i,1}F_1f_{4,j}
\]
where $(i,j)$ range over the set of indices for the block
containing $F_1$ and so
\[
\psi(a_1) = \sum a^1_{i,j}\psi(f_{i,1})\psi(F_1)\psi(f_{4,j})
\]
Since the matrix units $f_{i,1}, f_{4,j}$ lie in the selfadjoint
subalgebra $A_1 \cap A_1^*$, their images under $\psi$ are
sums of matrix units and the matrix above is identifiable,
after conjugation by a permutation unitary, with
\[
\left[ \begin{array}{ccc|ccc|ccc}
\alpha_1 \otimes a_1 & & 0 & & & & \lambda_1 \otimes a_1 &  & 0 \\
0 & 0 & & & & & 0 &  0 & \\
& 0& 0 & & & & & 0 & 0  \\ \hline
\beta_1 \otimes a_1 & & 0 & \gamma_1 \otimes a_1 & & 0 \\
0 & 0 & & 0 & 0 & \phantom{2} & \phantom{2} \\
& 0 & 0 & & 0 & 0\\ \hline
& & & \delta_1 \otimes a_1 & & 0 & \epsilon_1 \otimes a_1  & & 0 \\
& & & 0 & 0 &  & 0 & 0 \\
& & & & 0 & 0  & & 0 & 0
\end{array} \right]
\]

Thus, after conjugation by a permutation unitary, we can arrange
$\psi \circ \phi$ so that,
for each $i$, $\psi\circ\phi(E_i)$ has off-diagonal part
$$
\left[ \begin{array}{ccc|ccc|ccc}
\alpha_1 \otimes a_i & & \alpha_6 \otimes f_i & & & & \lambda_1 \otimes  
a_i &
	& \lambda_6 \otimes f_i \\
\alpha_2 \otimes b_i & \alpha_3 \otimes c_i & & & & & \lambda_2 \otimes b_i & 
	\lambda_3 \otimes c_i & \\
& \alpha_4 \otimes d_i & \alpha_5 \otimes e_i & & & & & \lambda_4
\otimes d_i &
	\lambda_5 \otimes e_i  \\ \hline
\beta_1 \otimes a_i & & \beta_6 \otimes f_i & \gamma_1 \otimes a_i & &
	\gamma_6 \otimes f_i \\
\beta_2 \otimes b_i & \beta_3 \otimes c_i & & \gamma_2 \otimes b_i &
\gamma_3 \otimes c_i & \phantom{2} & \phantom{2} \\
& \beta_4 \otimes d_i & \beta_5 \otimes e_i & & \gamma_4 \otimes d_i &  
\gamma_5 \otimes e_i \\ \hline
& & & \delta_1 \otimes a_i & & \delta_6 \otimes f_i & \epsilon_1 \otimes a_i 
	& & \epsilon_6 \otimes f_i \\
& & & \delta_2 \otimes b_i & \delta_3 \otimes c_i &
	& \epsilon_2 \otimes b_i & \epsilon_3 \otimes c_i \\
& & & & \delta_4 \otimes d_i & \delta_5 \otimes e_i
	& & \epsilon_4 \otimes d_i & \epsilon_5 \otimes e_i
\end{array} \right].
$$

As $(S \otimes T)^* = S^* \otimes T^*$, it follows that if $S_i$ and
$T_i$ are
partial isometries, $i=1,2$, and if $S_1$ and $S_2$ have orthogonal
initial projections
then so do $S_1 \otimes T_1$ and $S_2 \otimes T_2$.
Similar statements apply for final projections and for $T_1$ and $T_2$. 

In particular, consider the initial projection of $\alpha_1 \otimes a_1$. 
As $\psi \circ \phi$ is rigid, the (proper) 6-cycle $E_1, \ldots, E_6$ is 
sent to a direct sum of (proper) rank-one 6-cycles in $A_3$.
As $\alpha_1 \otimes a_1$ is a restriction of $\psi \circ \phi(E_1)$,
its initial projection is the initial projection of some restriction
of $\psi\circ\phi(E_2)$.
By the rigidity of $\psi \circ \phi$, this latter restriction must be some  
combination of $\beta_1 \otimes a_2, \ldots, \beta_6 \otimes f_2$.
However, the initial projection of $a_1$ is orthogonal to those of
$c_2$, $d_2$, $e_2$, and $f_2$, and
$\alpha_1$ has initial projection
orthogonal to that of $\beta_1$, so the only possibility is
that this latter restriction is a subprojection of $\beta_2 \otimes b_2$.
In fact, we must have equality, since we can argue reciprocally.

Since $\alpha_1 \otimes a_1$ and $\beta_2 \otimes b_2$ have coincidental 
initial projections it follows that $a_1$ and $b_2$ have coincidental initial
projections.
Continuing in this way for $\beta_2 \otimes b_2$ and $\gamma_3 \otimes c_3$,
we have that $a_1, b_2, c_3, d_4, e_5, f_6$ is a 6-cycle of
partial isometries in $A_2$.
This argument can be applied to each $a_i$ and so the $\phi(E_i)$ is
a direct sum of 6-cycles, showing $\phi$ is rigid.

Returning to $\alpha_1 \otimes a_1$ and $\beta_2 \otimes b_2$, we also have 
that $\alpha_1$ and $\beta_1$ have the same initial projection.
Repeating the argument again, it follows that $\psi$ is also rigid.
\end{proof}

The importance of the last two lemmas is that they immediately
give the following theorem.
In particular, invariants for regular isomorphisms of the
towers (with rigid embeddings) are in fact invariant for
the associated locally finite algebras.

\begin{thm}
\label{thm:algrigid}
Let $n \ge 3$.
Suppose that
\[A_1 \subseteq A_2 \subseteq \cdots \hbox{ and }
A'_1 \subseteq A'_2 \subseteq \cdots \]
are towers of direct sums of $2n$-cycle algebras,
where the inclusions are rigid.
Then the locally finite algebras $\sA_0 = \indalglim A_k$
and $\sA'_0 = \indalglim A_k'$ are star extendibly
isomorphic by a map $\phi$ with inverse $\psi=\phi^{-1}$
if and only there is commuting diagram
$$ \begin{CD}
A_1 @>>> A_{m_1} @>>> A_{m_2} @>>> A_{m_3} @>>> \cdots & \sA_0 \\
& \symbse{\phi_1} & \symbup{\psi_1} & \symbse{\phi_2} & \symbup{\psi_2}
      & \symbse{\phi_3} & \symbup{\psi_3} & \symbse{\phi_4}
      & & \symbupdown{{\phi,\psi}}\\
A'_1 @>>> A'_{n_1} @>>> A'_{n_2} @>>> A'_{n_3} @>>>  \cdots & \sA_0'
\end{CD} $$
where all the linking maps $\phi_k$ and $\psi_j$ are rigid and
$\phi = \indlim \phi_k$.
\end{thm}

\begin{proof}
One direction is clear and for the other,
assume that $\phi : \sA \to \sA'$ is given.
Then, since $A_1$ is finitely generated, there
exists $n_1$ and $m_1$ so that
\[ \phi(A_1) \subseteq A_{n_1}' \hbox{ and }
\psi(A_{n_1}') \subseteq A_{m_1}. \]
The restriction embeddings, $\phi_1 : A_1 \to A_{n_1}'$
and $\psi_1 : A_{n_1}' \to A_{m_1}$, are star
extendible and by hypothesis, have a rigid composition.
Now Lemmas~\ref{lma:PIs} and~\ref{lma:rigidfactors}
apply, to show that $\phi_1$ and $\psi_1$ are
rigid.
Continuing in this way, we obtain the required
diagram.
\end{proof}

\section{Approximate Factorisations}

To classify norm-closed direct limits, we need approximate versions of
Lemmas~\ref{lma:PIs} and~\ref{lma:rigidfactors}.
The latter will be used to define homology invariants for
operator algebra direct limits.
Readers interested only in the classification of algebraic direct limits
may go to Section 4.
For clarity and convenience, we have chosen not to determine the absolute
dependence of $\delta$ upon $\epsilon$ and $n$ in the proofs below.

\begin{lma}
\label{lma:approxPIs}
Suppose $\epsilon > 0$ and $a = \bigl( a_{ij} \bigr)$ is a block matrix
in a $2n$-cycle algebra $A$, $n \ge 3$.
There is constant $\delta >0$ so that if $a$ is $\delta$-close to
a partial isometry with initial and final projections in $A$, then
each entry $a_{ij}$ is $\epsilon$-close to a partial isometry.
\end{lma}

The proof is routine.

\begin{lma}
\label{lma:approxrigid}
Let $n \ge 3$, $A_1,A_2,A_3$ be $2n$-cycle algebras, and $\epsilon > 0$.
Let $\phi: A_1 \to A_2$ and $\psi : A_2 \to A_3$ be linear injections
which are $\delta$-close to star-extendible embeddings
and are such that
the composition $\psi \circ \phi$ is $\delta$-close to a rigid
embedding, $i : A_1 \to A_3$.
Then there is a sufficiently small $\delta$, depending on $\epsilon$ and $n$,
so that $\phi$ and $\psi$ are $\epsilon$-close to rigid embeddings.
\end{lma}

\begin{proof}
First observe that $A_1$ admits a direct sum decomposition
$A_1 = A_1 \cap A_1^* + \rad(A_1)$ and $\phi(a_1+a_2)=\phi_1(a_1)+\phi_2(a_2)$
where $\phi_1$ maps $A_1 \cap A_1^*$ to $A_2$ and $\phi_2$ maps $\rad(A_1)$
to $A_2$.
By standard finite-dimensional \cstar-algebra
theory~\cite[Chapter~III]{Davidson96}
$\phi_1$ is close to a star-extendible injection, so we may replace
it by a nearby \cstar-algebra homomorphism
$\phi_1' : A_1 \cap A_1^* \to A_2 \cap A_2^*$.
Similarly, we may assume the restriction of $\psi$ to $A_2\cap A_2^*$
is a \cstar-algebra homomorphism.

Assume now that the sets of matrix units $e^k_{ij}$, $1 \le k \le 3$,
have been chosen for $A_1$, $A_2$, and $A_3$ so that $\phi$ and $\psi$
map matrix units of $A_1 \cap A_1^*$ and $A_2 \cap A_2^*$ to sums of
matrix units.
Since $\phi$ is close to star-extendible injection, $\phi(e_{ij}^1)$
is close to $v_{ij}^1=\phi(e^1_{ii})\phi(e^1_{ij})\phi(e^1_{jj})$.
Let $\phi'$ be the linear map for which $\phi'(e^1_{ij})=v^1_{ij}$.
Then $\phi'$ agrees with $\phi$ on $A_1 \cap A_1^*$ and $\phi'$ is a
$C_1$-bimodule map, where $C_1$ is the diagonal masa spanned by $\{e_{ii}^1\}$
in $A_1$.
Define $\psi'$ similarly.
Since $\phi$ and $\psi$ are close to $\phi'$ and $\psi'$, we need only
prove the lemma for such maps.

Assume that $n=3$; this simplifies the notation.
The same argument applies for $n > 3$, with a smaller $\delta$.

We can now write, as in the proof of Lemma~\ref{lma:rigidfactors},
$$
\phi(E_i) = \left[ \begin{array}{ccc} a_i & & f_i \\ b_i & c_i & \\
	& d_i & e_i \end{array} \right], \qquad
\psi(F_i) = \left[ \begin{array}{ccc} \alpha_i & & \lambda_i \\ \beta_i  
	& \gamma_i & \\ & \delta_i & \epsilon_i \end{array} \right],
$$
for $i=1,2,\ldots,6$, where the $E_i$ and $F_i$ are rank-one six cycles
in $A_1$ and $A_2$ respectively.
These operators are close to partial isometries and, moreover, by
assumption, $\phi$ and $\psi$ are close to star-extendible homomorphisms.
By Lemma~\ref{lma:approxPIs}, each of the entries of each
$\phi(E_i)$ and $\psi(F_i)$ is close to a partial isometry.
Thus, we may assume that each of these entries is a  partial
isometry and that we still have $\phi$, $\psi$ close to
star-extendible homomorphisms.
Since $\beta_1\alpha_1^*$ has small norm (as $\psi(F_1)$ is close to
a partial isometry) and $\psi\circ\phi$ is close to a rigid
projection, it follows that the initial projection of $\alpha_1\otimes a_1$
is close to the initial projection of $\beta_2 \otimes b_2$ and hence that
the initial projection of $a_1$ is close to the initial projection of $b_2$.

Now redefine $\phi$ to obtain a new $\phi$ for which, in addition to
our earlier
assumptions, the initial and final projections of $a_1,b_2,c_3,d_4,e_5$ 
match up.
Also, redefine $f_6$ as $a_1b_2^*c_3d_4^*e_5$, which is close to the
original $f_6$ because $\phi$ is close to a star-extendible embedding.
The map $\phi$ now has a summand, determined by the $6$-cycle
$a_1,\ldots,f_6$,
which is a rigid algebra homomorphism.
Continuing in this way, we can construct a rigid algebra homomorphism
which is close to the original $\phi$.
Performing similar constructions with $\psi$, it follows that $\psi$ is
also close to a rigid embedding.
\end{proof}

\begin{thm}
\label{thm:approxalgrigid}
Let $n \ge 3$.
Suppose that
\[A_1 \subseteq A_2 \subseteq \cdots \hbox{ and }
A'_1 \subseteq A'_2 \subseteq \cdots \]
are towers of direct sums of $2n$-cycle algebras,
where the inclusions are rigid.
Then the operator algebras $\sA = \indlim A_k$
and $\sA' = \indlim A_k'$ are star extendibly
isomorphic if and only if there is commuting diagram
$$ \begin{CD}
A_1 @>>> A_{m_1} @>>> A_{m_2} @>>> A_{m_3} @>>> A_{m_4} @>>> \cdots & \sA \\
& \symbse{\phi_1} & \symbup{\psi_1} & \symbse{\phi_2} & \symbup{\psi_2}
      & \symbse{\phi_3} & \symbup{\psi_3} & \symbse{\phi_4} & \symbup{\psi_4}
      & \symbse{\phi_5} \\
A'_1 @>>> A'_{n_1} @>>> A'_{n_2} @>>> A'_{n_3} @>>> A'_{n_4} @>>>
\cdots & \sA'
\end{CD} $$
where all the linking maps $\phi_k$ and $\psi_j$ are rigid.
\end{thm}

\begin{proof}
Let $\Phi : \sA \to \sA'$ be a star-extendible isomorphism.
Then the restriction $\phi_0 = \Phi |_{A_1}$ is close to a
linear injection $\phi : A_1 \to A_k'$, for suitably large $k$.
Similarly, $\psi_0 = \Phi^{-1} |_{A_k'}$ is close to a linear injection
$\psi : A_k' \to A_l$ for suitably large $l$, and the composition,
$\psi \circ \phi$ is close to the given rigid embedding $A_1 \to A_l$.

We cannot immediately apply Lemma~\ref{lma:approxrigid} as stated,
since we only know that $\phi$ is to close to a star-extendible embedding
from $A_1$ to $\sA'$, and not one from $A_1$ to $A_k'$.
However, the same proof applies to this case, so we can apply the
lemma.
Thus, $\phi$ and $\psi$ are close to rigid embeddings, $\phi_1$ and
$\psi_1$ say.
Since $\psi_1 \circ \phi_1$ is close to the given rigid embedding,
it follows that they have the same  induced maps on $K_0$ and $H_1$
and so by Lemma 4.2 below are inner conjugate.
Adjusting $\phi_1$ by a unitary in $A_{n_1}$, we obtain a commuting triangle.
Continuing in the usual way, we can build the required diagram.
\end{proof}

The argument above shows more than we have stated, namely that each
star-extendible isomorphism $\Phi$ from $\sA$ to $\sA'$ determines a
star-extendible isomorphism $\{\phi_k,\psi_k\}$ for the towers,
which is unique up to an approximately inner automorphism, that is,
to a pointwise limit of unitary automorphisms.
Finally, if we define $H_1\sA$ to be the limit group $H_1\sA_0$
for the locally finite algebras $\sA_0$,
then we immediately have the following corollary.

\begin{cor}
The abelian group  $H_1\sA$ is well-defined and is
an invariant for star-extendible isomorphisms.
\end{cor}

\section{Invariants and Classification}


We now develop the invariants and prove the main theorem.

To each rigid embedding between $2m$-cycle algebras we may associate
an ordered $2m$-tuple, $(r_1,\ldots,r_{2m})$, which we call the
\textsl{multiplicity signature}, which is the set of multiplicities
of the $2n$ classes of multiplicity one embeddings in the direct
sum decomposition of $\phi$.
This signature depends on the identification of the reduced digraphs
of $A_1$ and $A_2$ with $D_{2m}$ and the labeling of the automorphisms
of $D_{2m}$.
Plainly, embeddings $\phi$, $\phi'$ are conjugate (inner unitarily
equivalent) if and only if they have the same multiplicity signature.

For the scaled $K_0$-group classification of limits of finite-dimensional
semisimple algebras, the cornerstone lemma is that two embeddings between
such algebras are conjugate if and only if they induce the same $K_0$-group
homomorphism.
Not surprisingly, the $K_0$-group homomorphism does not determine the
conjugacy class of embeddings between $2m$-cycle algebras.
However, if $H_1\phi$ denotes the group homomorphism $\bZ \to \bZ$
given by multiplication by
$$ r_1 - r_2 + r_3 - \cdots - r_{2m}$$
then, as we see below,  $K_0\phi$ and $H_1\phi$ together do determine the
multiplicity signature of $\phi$ and hence the conjugacy class of $\phi$.

It is natural then to seek a complete classification of the locally
finite algebras of Theorem~\ref{thm:algrigid} in terms of $K_0 A$
and the abelian group
$$ H_1 A := \indlim( \bZ, H_1\alpha_i) $$
That this group is indeed an invariant for star-extendible
isomorphism follows from Theorem~\ref{thm:algrigid}.

It has already been made clear, in the consideration of 4-cycle algebras
in Donsig and Power~\cite{DonPow97} and in Power~\cite{PowPRc},
that, beyond the scaled ordered group $K_0 A$ and the abelian group
$H_1 A$, it is necessary to consider a number of other invariants.
For some subfamilies of algebras, these additional invariants simplify or 
vanish but in the most general case, the appropriate invariant is a joint 
scale, $\Sigma A$ in $K_0 A \oplus H_1 A$, which embodies these additional
invariants.

\begin{defn}
Let $A$ be a locally finite algebra, as in Theorem~\ref{thm:algrigid}.
Then the \textsl{joint scale}, $\Sigma A$, is the subset of
$K_0 A \oplus H_1 A$ given by the elements of the form
$$ (K_0 \phi(e_{11}\oplus e_{22}), H_1\phi(g) )$$
where $\phi : A(D_{2m}) \to A$ is a rigid embedding and $g$ is a fixed
generator for $H_1(A(D_{2m})) = \bZ$.
Furthermore, we define the unital joint scale
to be the subset arising from unital embeddings $\phi$.
\end{defn}

It is possible to give an intrinsic formulation of $H_1 A$ for a digraph
algebra $A$ and hence an intrinsic formulation of $H_1 \phi$, although
we do not do so here.

\begin{lma}
If $A$ is $2m$-cycle algebra and $\phi, \psi$ are rigid embeddings
from $A(D_{2m})$ to $A$ with
$$  (K_0 \phi(e_{11}\oplus e_{22}), H_1\phi(g) )
	= (K_0 \psi(e_{11}\oplus e_{22}), H_1\psi(g) ), $$
then $\phi$ and $\psi$ are conjugate.
\end{lma}

This lemma has essentially the same proof as the analogous 4-cycle
result \cite[Lemma~11.4]{DonPow97}; see also
\cite[Proposition~3.4]{Pow96a}.

It is instructive to note the possible variation in $H_1 \phi$ once
$K_0 \phi$ is specified.
If $\phi$ has multiplicity signature $(r_1,\ldots,r_{2m})$,
and $K_0\phi = K_0\psi$, then for suitable integers $k$,
$\psi$ has multiplicity signature
$$ (r_1+k,r_2-k,r_3+k,\ldots, r_{2m}-k)$$
and $H_1\psi = H_1\phi + [2mk]$.
The range for $k$ is $-M_2 \le k \le M_1$ where $M_1$ and $M_2$
are the maximum and minimum values determining nonnegative signature.
The set
$$ \hr(\phi) = \{ H_1 \psi : \psi \text{ rigid and } K_0\psi = K_0 \phi \}$$
is called the \textsl{homology range} of $\phi$ and is realized as an interval
in the mod $2m$ congruence class $M_1 + 2m\bZ$.
These remarks already suggest that there can be congruence class obstructions
to the lifting of $K_0H_1$ group isomorphisms
$$ \gamma_0 \oplus \gamma_1 : K_0 A \oplus H_1 A \to K_0 A' \oplus H_1 A' $$
to algebra isomorphisms. (For further discussions in the case of
4-cycle algebras see \cite{DonPow97}.)

\begin{defn}
Let $n \ge 3$ and suppose that
\[A_1 \subseteq A_2 \subseteq \cdots \hbox{ and }
A'_1 \subseteq A'_2 \subseteq \cdots \]
are towers of direct sums of $2n$-cycle algebras,
where the inclusions are rigid.
Letting $\sA = \indlim A_k$ and $\sA' = \indlim A_k'$,
we say that a scaled, ordered group isomorphism
\[ \gamma_0 : K_0 \sA \to K_0 \sA' \]
is \textit{of rigid type} if there is some
rigid isomorphism $\Gamma : \sA \to \sA'$ so that
$K_0 \Gamma = \gamma_0$. Equivalently $\gamma_0$ is induced by
a commuting diagram of
rigid embeddings (as in Theorem 2.4).
\end{defn}

At the finite-dimensional level, a $K_0$-group homomorphism, of
$6$-cycle algebras is of rigid type if it can be realized as a sum of two
integral matrices of the form:
\[ \begin{bmatrix} a & b & c & 0 & 0 & 0 \\  c & a & b & 0 & 0 & 0 \\
	b & c & a & 0 & 0 & 0 \\ 0 & 0 & 0 & a & b & c \\
	0 & 0 & 0 & c & a & b \\ 0 & 0 & 0 & b & c & a \end{bmatrix}
   + \begin{bmatrix} d & e & f & 0 & 0 & 0 \\  e & f & d & 0 & 0 & 0 \\
	f & d & e & 0 & 0 & 0 \\ 0 & 0 & 0 & e & f & d \\
	0 & 0 & 0 & f & d & e \\ 0 & 0 & 0 & d & e & f \end{bmatrix} \]
Similar sums describe the inclusions for larger $2n$-cycle algebras.

We can now obtain the main theorem of the paper, which
(but for the case of even 4-cycle systems)
completes the classification of operator algebras begun
in Power~\cite[Chapter~11]{Power92} and continued in
Donsig and Power~\cite{DonPow97} and in Power~\cite{PowPRc}.

\begin{thm}
\label{thm:mainclass}
Let $n \ge 3$ and suppose that
\[A_1 \subseteq A_2 \subseteq \cdots \hbox{ and }
A'_1 \subseteq A'_2 \subseteq \cdots \]
are towers of direct sums of $2n$-cycle algebras,
where the inclusions are rigid.
Then the following conditions are equivalent:
\begin{enumerate}
\item There is a commuting diagram
$$ \begin{CD}
A_1 @>>> A_{m_1} @>>> A_{m_2} @>>> A_{m_3} @>>> \cdots & \sA \\
& \symbse{\phi_1} & \symbup{\psi_1} & \symbse{\phi_2} & \symbup{\psi_2}
      & \symbse{\phi_3} & \symbup{\psi_3} & \symbse{\phi_4}
      & & \symbupdown{{\phi,\psi}}\\
A'_1 @>>> A'_{n_1} @>>> A'_{n_2} @>>> A'_{n_3} @>>> \cdots & \sA'
\end{CD} $$
where all the linking maps $\phi_i$ and $\psi_j$ are rigid.
\item The locally finite algebras $\sA_0 = \indalglim A_k$
and $\sA_0' = \indalglim A_k'$ are star-extendibly isomorphic.
\item The operator algebras $\sA = \indlim A_k$
and $\sA' = \indlim A_k'$ are star extendibly isomorphic.
\item There is an abelian group isomorphism
$$ \gamma_0 \oplus \gamma_1 :
	K_0 \sA \oplus H_1 \sA \to K_0 \sA' \oplus H_1 \sA', $$
where $\gamma_0$ is a scaled, ordered group isomorphism of rigid type
and $\gamma_0 \oplus \gamma_1$ preserves the joint scale. Also, in
the unital case it suffices to preserve the unital joint scale.
\end{enumerate}
\end{thm}

\begin{proof}
 From Theorems~\ref{thm:algrigid} and~\ref{thm:approxalgrigid},
we know that the first three conditions are equivalent.
Clearly, the first implies the fourth, so it remains only to
show that the fourth implies the first.
The argument follows a familar scheme, seen for example
in~\cite[Theorem~11.5]{DonPow97} and~\cite[Theorem 5.2]{PowPRc},
but for completeness, we outline the argument.

Consider a multiplicity one
rigid embedding $\phi : A(D_{2n}) \to A_1$, which determines
an element $([\phi(e_{11}\oplus e_{22})],\delta)$
of the joint scale of $A_1$.
By condition~4, there is some rigid embedding $\eta : A(D_{2n}) \to A_k'$,
for some $k$, so that
$$\gamma_0 \oplus \gamma_1 ([\phi(e_{11}\oplus e_{22})],\delta) =
	([\eta(e_{11}\oplus e_{22})], H_1\eta(g)). $$
We may assume that matrix units have been chosen for $A_k'$ so that
$\eta$ maps matrix units to sums of matrix units.

We claim that $\eta$ has an extension $A_1 \to A_k'$ (after possibly
increasing $k$), which we also denote $\eta$.
To see this first observe that since $\gamma_0$ is a scaled group
isomorphism there is a C*-algebra homomorphism
$$ \eta_0 : A_1 \cap A_1^* \to A_k' \cap (A_k')^* $$
so that $K_0\eta_0$ is equal to $\gamma_0$ restricted to $K_0 A_1$
and, moreover, $\eta_0$ maps matrix units to sums of matrix units.
Since each matrix unit $e \in A_1$ admits a unique factorization
$e=e_1fe_2$ where $e_1,e_2 \in A_1 \cap A_1^*$ and $f$ is in
$\phi(A(D_{2n}))$
it follows that there is a unique star extendible extension of $\eta$
and $\eta_0$ to a rigid embedding $\zeta : A_1 \to A_k'$.
Although $\eta$ depends on $\phi$ and $\eta_0$ the inner conjugacy class
of $\eta$ is determined.


Next, consider a multiplicity-one rigid embedding $\psi: A(D_{2n}) \to A_k'$
which induces an element of the joint scale
$([\psi(e_{11}\oplus e_{22})], H_1\psi(g))$.
We assume for convenience that matrix units have been chosen
for $A_k'$ so that $\psi$ maps matrix units to sums of matrix units.
Again, there is some rigid embedding $\zeta : \psi(A(D_{2n})) \to A_l$,
for some $l$, so that
$$\gamma_0^{-1} \oplus \gamma_1^{-1} ([\psi(e_{11}\oplus e_{22})],
	H_1\psi(g)) = ([\zeta(e_{11}\oplus e_{22})], H_1\zeta(g)). $$
We may assume that matrix units have been chosen for $A_l$ so that
$\zeta$ maps matrix units to sums of matrix units.

As with $\eta$, $\zeta$ has a star extendible extension $A_k' \to A_l$
(after possibly increasing $l$), which we also denote $\zeta$.
Since $\zeta \circ \eta$ and $\alpha_{l-1} \circ \cdots \circ \alpha_1$
induce the same $K_0\oplus H_1$ map, we may replace $\zeta$ by an
inner conjugate map so that the algebra maps are equal.
Continuing in this way, we can construct the required diagram.
\end{proof}

For an elementary illustration of the theorem we now consider unital
limit algebras for which $A\cap A^*$ is a direct sum of two UHF C*-algebras.
We say such algebras are \textit{of matroid type}.
For a fixed positive integer $d$ consider
the stationary systems of the form
$$
\begin{CD}
A(D_6) @>>> A(D_6)\otimes M_{3d} @>>> A(D_6)\otimes M_{(3d)^2} @>>>   
\cdots  \\
\end{CD}
$$
where the $k$-th embedding $\alpha_k$ satisfies
\[
K_0(\alpha_k) = d
\begin{bmatrix} 1 & 1 & 1 & 0 & 0 & 0 \\  1 & 1 & 1 & 0 & 0 & 0 \\
        1 & 1 & 1 & 0 & 0 & 0 \\ 0 & 0 & 0 & 1 & 1 & 1 \\
        0 & 0 & 0 & 1 & 1 & 1 \\ 0 & 0 & 0 & 1 & 1 & 1 \end{bmatrix},
~H_1(\alpha_k) = [s].
\]
There are $d + 1$ possible values of $s$,
namely the integers that lie in the homology range
\[
S = \{-3d, -3d+6, \dots, 3d -6, 3d\}
\]
and so there are $d+1$ associated operator algebras
$A_s, s \in S$.

The group $K_0(A_s)$ identifies naturally with
the subgroup of $\bQ \oplus \bQ$ corresponding to the generalised integer
${3^\infty d^\infty}$
with order unit $1 \oplus 1$ whilst $H_1(A_s)$ identifies
naturally with the subgroup determined by $s^\infty$.
In the exceptional case $s=0$, $H_1(A_s)=0$ and
so the joint scale then reduces to the usual scale in $K_0(A_s)$.

In general therefore, the unital part of the joint scale is the subset of 
\[
{\bZ}[\frac{1}{3^\infty d^\infty}] \oplus
{\bZ}[\frac{1}{3^\infty d^\infty}] \oplus {\bZ}[\frac{1}{s^\infty}]
\]
arising from unital rigid embeddings
$\phi : A(D_6) \to A_s$.
In the case of the extreme values of $s$ the
algebras $A_{3d}, A_{-3d}$,
are isomorphic and the unital joint scale coincides with
the subset of elements $ 1/3 \oplus 1/3 \oplus h$
where for some positive integer $m$ one has
$h = k/(3d)^m$ where $k \equiv (3d)^m $ mod $6$
and $-(3d)^m \le k \le (3d)^m$.
In particular the homology component of the
unital joint scale is a symmetric set in a
symmetric finite interval of $H_1A_s$.
In contrast, in the nonextreme case, $s \ne 3d, -3d$,
the unital joint scale is identified with the set of elements
for which $h$ merely satisfies the congruence restriction.
(In our simple example the congruence restriction has no
consequence if $d$ is even whilst if $d$ is odd then
the homology part of the unital joint scale corresponds to
odd numerators.) Indeed,
for such an element $h$ in ${\bZ}[\frac{1}{s^\infty}]$ we may
choose $t$ large
enough so that
\[
h = \frac{k}{(3s)^m} = \frac{k'}{(3s)^t}
\]
where $ k' \le (3d)^t$. Then there is a unital embedding
\[
\psi : A(D_6) \to A(D_6) \otimes M_{(3d)^t}\]
with $H(\psi) = [k']$, and so $\psi$
determines the element $1/3 \oplus 1/3 \oplus h$.

The $H_1$ group is already a distinguishing
invariant for $A_{s_1}, A_{s_2}$ if $s_1$ and $s_2$ have different prime
divisors. On the other hand if $s_1$ and $s_2$ are nonextreme
and have the same prime divisors then one can verify
that the unital joint scales
coincide and so, by the theorem, $A_{s_1}, A_{s_2}$ are isomorphic.

It follows from similar observations that if $A, A'$
are unital stationary limits of $2n$-cycle algebras,
of nonextreme type,
then $A, A'$ are star extendibly isomorphic if and only if
$C^*(A)$ and $C^*(A')$
are isomorphic and $H_1(A)$ and $H_1(A')$ are isomorphic
abelian groups.
In fact the same conclusion is true for the more general limits of
``homologically limited" systems of matroid type, that is, those for which 
the natural scale of $H_1$ is not a finite interval but, as above,
coincides with $H_1(A)$.

\providecommand{\bysame}{\leavevmode\hbox to3em{\hrulefill}\thinspace}

\end{document}